\documentclass[12pt]{amsart}

\newsymbol\pp 1275
\newsymbol\twoheadrightarrow 1310

\newcommand{\letters}{\renewcommand{\theenumi}{\alph{enumi}}}
\newcommand{\denom}{\operatorname{denom}}
\newcommand{\id}{\operatorname{id}}
\newcommand{\udenom}{\operatorname{udenom}} 
\newtheorem{theorem}{Theorem}[section]
\newtheorem{proposition}[theorem]{Proposition}

\newtheorem{lemma}[theorem]{Lemma}
\newtheorem{corollary}[theorem]{Corollary}

\theoremstyle{definition}

\newtheorem{definition}[theorem]{Definition}
\newtheorem{remark}[theorem]{Remark}
\newtheorem{example}[theorem]{Example}

\newcommand{\Q}{{\mathbb Q}}

\newcommand{\cx}{{\mathcal C}}
\newcommand{\Z}{{\mathbb Z}}

\newcommand{\SL}{\operatorname{SL}}

\newcommand{\R}{{\mathbb R}}
\newcommand{\rank}{\operatorname{rank}}

\newcommand{\lcm}{\operatorname{lcm}}

\newcommand{\N}{{\mathbb N}}

\title{Universal Denominators of Hilbert Series}

\author{Harm Derksen}
\thanks{The author is supported by NSF grant 0102193.}
\begin{document}

\maketitle
\begin{abstract}

The denominator
of the Hilbert series of a finitely generated $R$-module $M$ 
does not always divide the denominator of the
Hilbert series of $R$. For this reason, we define the universal denominator. 
The universal denominator of a module $M$ is
  the least common multiple
of the denominators of the Hilbert series of all submodules of $M$.
The universal denominator behaves nicely with respect to short
exact sequences and tensor products. It also has interesting geometric 
interpretations.
Formulas are given for the universal denominator for rings
of invariants. Dixmier gave a conjectural formula
for the denominator of the Hilbert series
of invariants of binary forms. We show that
the universal denominator 
is actually equal to Dixmier's formula in that case.
\end{abstract}
\section{Definitions and basic properties}
Let $\N=\{0,1,2,\dots\}$ be the set of natural numbers and let $K$
be the  base field. With a multi-graded ring of finite type we mean
a multi-graded ring
 $R=\bigoplus_{d\in \N^r}R_d$ such that $R_0$ is a finite dimensional 
 $K$-vector space and $R$ is finitely generated over $K$.
\begin{definition}
If $M=\bigoplus_{d\in \Z^r}M_d$ is a finitely generated multi-graded 
$R$-module then $M_d$ is finite dimensional for all $d$. We define 
the Hilbert series of $M$ by
$$
H(M,t)=\sum_{d\in \Z^r}\dim(M_d)\,t^d
$$
where $t=(t_1,\dots,t_r)$, $d=(d_1,\dots,d_r)$ and
$t^d=t_1^{d_1}t_2^{d_2}\cdots t_r^{d_r}$.
\end{definition}
The Hilbert series $H(M,t)$ is a Laurent series in $t_1,t_2,\dots,t_r$.
The Hilbert series $H(R,t)$ of $R$ itself is a power series in $t_1,\dots,t_r$.
\begin{definition}
We can uniquely write $H(M,t)=A(t)/B(t)$ where
$A(t)$ is a Laurent polynomial in $t_1,\dots,t_r$,
$B(t)$ is a polynomial in $t_1,\dots,t_r$ 
with $B(0)=1$ and $A(t)$ and $B(t)$
do not have a common non-constant polynomial factor. We call $B(t)$ {\bf the} denominator
of $H(M,t)$, and we will denote it by $\denom(M,t)$.
\end{definition}
\begin{example}\label{ex1}
Let $R=K[x_1,\dots,x_n]$ be the multi-graded polynomial ring
where the variable $x_i$ has multidegree $d_i=(d_{i,1},\dots,d_{i,r})\in \N^r$
for $i=1,2,\dots, n$. The Hilbert series $H(R,t)$ is equal
to 
$$
\frac{1}{(1-t^{d_1})\cdots (1-t^{d_n})}.
$$
Hilbert proved that every finitely generated (multi-graded) $R$-module $M$
 has a finite resolution of finitely generated 
free $R$-modules (see~\cite{Hi}).His arguments
also showed that $H(M,t)$ of $M$ can be written as
$$
\frac{A(t)}{(1-t^{d_1})\cdots (1-t^{d_n})}.
$$
where $A(t)$ is a Laurent polynomial in $t_1,\dots,t_r$.
The denominator of $H(M,t)$ must divide
$$
(1-t^{d_1})\cdots (1-t^{d_n}).
$$
\end{example}
The argument in the previous example heavily uses the fact
that finitely generated modules over  graded polynomial rings
have finite free resolutions. The graded polynomial rings
are the only graded rings over $K$ for which the module $K$
has a finite free resolution (see~\cite{Serre}). For an arbitrary multi-graded ring $R$
of finite type
there is no reason why the denominator of the Hilbert series
of a finitely generated graded module should divide the
denominator of the Hilbert series of $R$.
This is indeed not always the case as the following example shows.
\begin{example}
Hilbert series of submodules and quotient modules may have a 
larger denominator. For example, take the
graded ring $R=K[x,y]/(y^2)$
where $x$ has degree 2 and $y$ has degree 1.
We have
$$
H(R,t)=\frac{1+t}{(1-t^2)}=\frac{1}{(1-t)},\quad \denom(R,t)=1-t
$$
Now $(y)$ is an ideal of $R$ and $R/(y)\cong K[x]$.
We get
$$
H((y),t)=\frac{t}{1-t^2},\quad\denom((y),t)=1-t^2
$$
and
$$
H(K[x],t)=\frac{1}{1-t^2},\quad\denom(K[x],t)=1-t^2.
$$
So the submodule $(y)$ and the quotient module $K[x]$ both
have Hilbert series with larger denominator. Also note that
the modules $(y)$ and $K[x]$ cannot have  finite free resolutions.
\end{example}

\begin{definition}
Let $B(t)$ be the unique smallest polynomial in $t_1,\dots,t_r$
 such that $B(0)=1$, 
and $B(t)H(N,t)$
is a Laurent polynomial for every multi-graded submodule $N\subseteq M$.
We call $B(t)$ the {\bf universal denominator} of $H(M,t)$,
and we denote it by $\udenom(M,t)$ (or $\udenom_R(M,t)$
since it may depend on $R$ if $M$ is a finitely generated
module for several choices of $R$).
\end{definition}
\begin{remark}\label{rem1}
From Example~\ref{ex1} follows that the universal denominator
is well defined, i.e., there exists
a polynomial $B(t)$ such that $B(t)H(N,t)$ is a Laurent
polynomial for all submodules $N$ of $M$. Indeed, if $R$ is generated by
homogeneous $f_1,\dots,f_s$ of degrees $d_1,\dots,d_s\in \N^r$ respectively,
then any finitely generated module of $R$ can be viewed as a 
finitely generated module
of $K[x_1,\dots,x_s]$ by the surjective ring homomorphism
$$
K[x_1,\dots,x_s]\to R=K[f_1,\dots,f_s],
$$
where $x_i\mapsto f_i$ for all $i$.
Every submodule $N$ of $M$ is a finitely generated
$K[x_1,\dots,x_s]$-module, and therefore
$$
(1-t^{d_1})\cdots (1-t^{d_s})H(N,t).
$$
is a Laurent polynomial. 
\end{remark}
\begin{lemma}
If 
$$
0\to M'\to M\stackrel{p}{\rightarrow} M''\to 0
$$
is an exact sequence of finitely generated graded $R$-modules
then
$$
\udenom(M,t)=\lcm(\udenom(M',t),\udenom(M'',t))
$$
where $\lcm$ is the least common multiple.
\end{lemma}
\begin{proof}
It follows easily from the definition that
$\udenom(M',t)$ and $\udenom(M'',t)$ divide $\udenom(M,t)$.
If $Z$ is a multi-graded submodule of $M$ then we have an exact sequence
$$
0\to M'\cap Z\to Z\stackrel{p}{\rightarrow} p(Z)\to 0.
$$
Since $H(Z,t)=H(M'\cap Z,t)+H(p(Z),t)$ it follows
that $\denom(Z,t)$ divides 
$$\lcm(\denom(M'\cap Z,t),\denom(p(Z),t)).$$
Since $M'\cap Z\subseteq M'$ and $p(Z)\subseteq M''$ for all
submodules $Z\subseteq M$,
we have that $\udenom(M,t)$ divides $\lcm(\udenom(M',t),\udenom(M'',t))$.
\end{proof}

\begin{definition}
For $d\in \N$, let $\phi_d(t)\in \Z[t]$ be the unique
(irreducible) polynomial with $\phi_d(0)=1$ and
whose zeroes are exactly the
primitive $d$-th roots of unity (up to 
a scalar $\pm 1$ this is the usual $d$-th cyclotomic polynomial). 
If $d=(d_1,\dots,d_r)\in \N^r$ and
$t=(t_1,\dots,t_r)$, then we define $\phi_d(t)=\phi_{k}(t^{d/k})$
where $k=\gcd(d_1,d_2,\dots,d_r)\in \N$ is the greatest common divisor.
\end{definition}
If $d,e\in \N^r$
then we say that $d$ divides $e$ if there exists a $k\in \N$ such that
$e=kd$. The least common multiple $\lcm(d,e)$ is
the smallest {\it nonzero\/} vector $f$ divisible
by $d$ and $e$
if such a vector $f$ exists. Otherwise $\lcm(d,e)$ is defined
to be the zero vector.
For example, we have $\lcm((4,2),(6,3))=(12,6)$ and 
$\lcm((4,2),(2,2))=(0,0)$.
\begin{definition}
Suppose now that $R$ is a multi-graded ring of finite type
over a field $K$. Let $I^{[d]}$ be the ideal of
$R$ generated by all $R_e$ for which $d$ does not divide $e$.
\end{definition}
\begin{theorem}\label{theo10}
Suppose that $M$ is a finitely generated graded $R$-module.
\begin{enumerate}
\letters
\item We have
$$
\udenom(M,t)=\prod_{d\in \N^r}\phi_d(t)^{m_d}
$$
where $m_d$ is the dimension of the support of $M/I^{[d]}M$.
\item 
The universal denominator $\udenom(M,t)$ is the greatest common
divisor of all
$$
(1-t^{d_1})(1-t^{d_2})\cdots (1-t^{d_s})
$$
for which there exist homogeneous $f_1,f_2,\dots,f_s\in R$
of degrees $d_1,d_2,\dots,d_s$ respectively such that $M$ is a finitely generated
$K[f_1,\dots,f_s]$-module.
\end{enumerate}
\end{theorem}
\begin{proof}
Define
$$
A(t)=\prod_{d\in \N^r}\phi_d(t)^{m_d}
$$
and let $B(t)$ be the greatest common divisor of all
$$
(1-t^{d_1})(1-t^{d_2})\cdots(1-t^{d_s})
$$
for which there exist homogeneous $f_1,\dots,f_s\in R$ of
degrees $d_1,d_2,\dots,d_s$ respectively such that
$M$ is a finitely generated $K[f_1,\dots,f_s]$-module. (We normalize
such that $B(0)=1$.)

The ring $R^{[d]}=R/I^{[d]}$ is graded, $R^{[d]}=\bigoplus_{k\geq 0}
R^{[d]}_{kd}$.
Let us define 
$$
M^{[d]}_{\overline{e}}:=\bigoplus_{a-e\in \Z d} M^{[d]}_a
$$
Now $M^{[d]}_{\overline{e}}$ is a $R$-submodule of $M^{[d]}$.
There is a finite decomposition
$$
M^{[d]}=\bigoplus_{\overline{e}\in \Z^r/\Z d}
M^{[d]}_{\overline{e}}.
$$
Because $M^{[d]}$ has $m_d$-dimensional support, $M^{[d]}_{\overline{e}}$
has $m_d$-dimensional support for some $e\in \Z^r$.
The Hilbert series $H(M^{[d]}_{\overline{e}},t)$ has the form
$t^eP(t^d)$ for some rational function $P(t)$. Therefore
$(1-t^d)^{m_d}$ divides $\denom(M^{[d]}_{\overline{e}},t)$,
$\udenom(M^{[d]},t)$ and $\udenom(M,t)$.
In particular, $\phi_{d}(t)^{m_d}$ divides $\udenom(M,t)$.
We conclude that $A(t)$ divides $\udenom(M,t)$.

From Remark~\ref{rem1} follows that
$\udenom(M,t)$ divides
$$(1-t^{d_1})(1-t^{d_2})\cdots (1-t^{d_s})$$
whenever there exist $f_1,f_2,\dots,f_s\in R$ as in the theorem.
This shows that $\udenom(M,t)$ divides $B(t)$.

Let us write 
$$
B(t)=\prod_{d\in \N^r} \phi_d(t)^{a_d}
$$
Choose $g_1,\dots,g_s\in I^{[d]}$ homogeneous generators.
Let $d_i$ be the degree of $g_i$. We may assume that $d$ does not
divide $d_i$ for all $i$.
Since $M^{[d]}$ has $m_d$-dimensional support, we may choose
$f_1,\dots,f_{m_d}$ homogeneous of degrees $e_1,\dots,e_{m_d}$
(all nonzero multiples of~$d$)
such that 
$$M^{[d]}/(f_1,\dots,f_{m_d})M^{[d]}=
M/(g_1,\dots,g_s,f_1\dots,f_{m_d})M$$ 
is finite dimensional. It follows that $M$ is a finite
$K[g_1,\dots,g_s,f_1,\dots,f_{m_d}]$-module by the homogeneous
Nakayama Lemma.
Now $B(t)$ divides 
$$(1-t^{d_1})\cdots (1-t^{d_s})(1-t^{e_1})\cdots (1-t^{e_{m_d}}).
$$
It follows that $a_d\leq m_d$. Since this holds for all $d$,
we get that $B(t)$ divides $A(t)$.

We have proven that $A(t)$ divides $\udenom(M,t)$, $\udenom(M,t)$ divides
$B(t)$ and $B(t)$ divides $A(t)$. Moreover, we have $A(0)=\udenom(M,0)=B(0)=1$.
It follows that
$$
\udenom(M,t)=A(t)=B(t).
$$
\end{proof}

\begin{definition}
Suppose that $R=\bigoplus_{d\in \N^r}R_d$
and $S=\bigoplus_{d\in \N^r}S_d$ are both finitely generated multi-graded rings.
Then we can define the graded tensor product 
$$R\otimes S=\bigoplus_{d\in
\N^r}(R\otimes S)_d$$
where 
$$
(R\otimes S)_d=\bigoplus_{e} R_e\otimes S_{d-e}.
$$
In a similar fashion we can define the tensor product of
a graded $R$-module with a graded $S$-module.
\end{definition}
\begin{lemma}
If $R$ and $S$ are as in the previous definition and $M$ and $N$
are finitely generated modules for $R$ and $S$ respectively, then
$$
\udenom_{R\otimes S}(M\otimes N,t)=\udenom_R(M,t)\udenom_S(N,t).
$$
\end{lemma}
\begin{proof}
Note that 
$$
(M\otimes N)^{[d]}=M^{[d]}\otimes N^{[d]}.
$$
Now it follows from Theorem~\ref{theo10}(a).
\end{proof}
\section{Geometry of the universal denominator}
Let us now give a geometric description 
of the universal denominator. For convenience, we will assume
that the base field $K$ is algebraically closed from now on.
Let $X$ be the
affine variety corresponding to $R$.
The multigrading on $R$ corresponds to the
action of an $r$-dimensional torus $T$ on $X$.
\begin{definition}
To each $d\in \N^r$ we can associate the character of $T$
defined by $t\mapsto t^d$.
Let $T^{[d]}$ be the kernel of this character and 
let $X^{[d]}$ be the zero set of the ideal $I^{[d]}$.
\end{definition}
\begin{lemma}
Suppose that $d=(d_1,\dots,d_r)$
and the characteristic of $K$ does not divide $d_i$ for some $i$.
Then the  set $X^{[d]}$ is equal to $X^{T^{[d]}}$, the fixed
point set of $T^{[d]}$.
\end{lemma}
\begin{proof}
Suppose that $x\in X^{[d]}$. Let $f\in R$ be homogeneous of degree $e$.
If $e$ is not divisible by $d$ then $f\in I^{[d]}$, so 
$f(x)=f(t\cdot x)=0$.
If $e$ is divisible by $d$, then
$$
f(x)=t^ef(x)=(t\cdot f)(x)=f(t^{-1}\cdot x)
$$
for all $t\in T^{[d]}$.
It follows that $f(x)=f(t\cdot x)$ for all $f\in R$, hence
$x=t\cdot x$ for all $t\in T^{[d]}$. This shows that $x\in X^{T^{[d]}}$.

Conversely, if $x\in X^{T^{[d]}}$, $f\in I^{[d]}$ is homogeneous
of degree $e$, and $e$ is not divisible by $d$,
then
$$
f(x)=f(t^{-1}\cdot x)=(t\cdot f)(x)=t^ef(x)
$$
for all $t\in T^{[d]}$. We can choose $t\in T^{[d]}$ such
that $t^e\neq 1$ (because $d$ does not divide $e$
and $d$ is not divisible by the characteristic of $K$).
It follows that $f(x)=0$.
We have that
$x\in X^{[d]}$, because $I^{[d]}$ is generated by such $f$.
\end{proof}
\section{the universal denominator in Invariant Theory}
The notion of the universal denominator is very useful
for Hilbert series of invariant rings. In fact, many of the
results in this paper where inspired by some arguments
in the paper \cite{DW} where the Hilbert series of invariant
rings of quiver representations were studied.

Suppose that $S=\bigoplus_{d\in \N^r}S_d$ is a 
multi-graded ring of finite type over $S_0=K$ and
suppose that a reductive linear algebraic group $G$ acts regularly on $S$
such that the action respects the multi-grading.
We know that the invariant ring $R:=S^G$ is finitely generated over $K$ as
well. As before we define $I^{[d]}\subset R$ to be the ideal
generated by all $R_e$ such that $d$ does not divide $e$.
Also define $J^{[d]}=SI^{[d]}$ be the ideal of $S$ generated by $I^{[d]}$.
Geometrically, let $X$ be the affine variety corresponding to $R$
and let $Y$ be the affine variety corresponding to $S$. Let $\pi:Y\to X$
be the categorical quotient map corresponding to the inclusion
$R=S^G\subseteq S$. The multigrading on $R$ and $S$ correspond to
the action of a torus $T$ on $X$ and on $Y$. The quotient map $\pi:Y\to X$
is $T$-equivariant. The zero set of $I^{[d]}$ is equal to $X^{[d]}$
which is the set of $T^{[d]}$ fixed points on $X$. The zero set of
$J^{[d]}=SI^{[d]}$ is $\pi^{-1}(X^{[d]})\subseteq Y$.
\begin{definition}
Let $\zeta\in T$ and let $g\in G$. Then we define
$$
Y^g_\zeta=\{y\in Y\mid g\cdot y=\zeta\cdot y\}.
$$
\end{definition}
\begin{remark}
 For some of the arguments that
 we are going to present it is useful to have
  an element $\zeta\in T$
such that $\zeta$ generates a dense subgroup of the
torus $T$. It is not
always possible to choose such an element. For example
if the base field is the algebraic closure of a finite field
then every element of $T$ will have finite order. On the other hand,
if $K$ contains $\Q$ or $K$ contains an element that is transcendent over
the prime field, then there will exists such an element
$\zeta\in T$ constructed as follows. Choose distinct normalized valuations
$v_1,v_2,\dots,v_r$. Choose $\zeta_i\in K^\star$
such that $v_i(\zeta_i)=1$ and $v_j(\zeta_i)=0$ for $j\neq i$.
Then the group generated by $\zeta=(\zeta_1,\dots,\zeta_r)\in T$ will lie dense in $T$.
By extending the base field with transcendental elements we may
always assume that there exists an element $\zeta\in T$ which
generates a dense subgroup. If $d\in \N^r$ is not divisible by 
the characteristic of $K$, then similar arguments show that $\zeta\in T^{[d]}$
can be chosen such that $\zeta$ generates a dense orbit of $T^{[d]}$.
\end{remark}
\begin{proposition}\label{prop101}
Suppose that $d$ is not divisible by
the characteristic of $K$ and suppose that  $\zeta\in T^{[d]}$
generates a dense subgroup of $T^{[d]}$.
If the orbit $G\cdot y$ is closed, then $\pi(y)\in X^{[d]}$ if
and only if $y\in Y^g_\zeta$ for some $g\in G$.
\end{proposition}
\begin{proof}
Suppose that $G\cdot y$ is closed. Then $\pi(y)\in X^{[d]}$ if and only
if 
$$\pi(\zeta\cdot y)=\zeta\cdot \pi(y)=\pi(y)$$
Every fiber of $\pi$ has only one closed orbit and $y$ and $\zeta\cdot y$
have closed orbits. In particular, $\pi(\zeta\cdot y)=\pi(y)$
if and only if $\zeta\cdot y=g\cdot y$ for some $g\in G$.
\end{proof}
\begin{theorem}\label{theo100}
Suppose that $G$ is a connected linearly reductive algebraic group and
$H$ is linearly reductive subgroup of $G$ containing a maximal torus of $G$. Then
$\udenom(S^G,t)$ divides $\udenom(S^H,t)$.
\end{theorem}
\begin{proof}
Let $Y$ be the variety corresponding to $S$
and let $\pi_G:Y\to Y/\!\!/G$ be the categorical quotient
with respect to $G$ (so $Y/\!\!/G$ is the variety corresponding
to $S^G$) and let  $\pi_H:Y\to Y/\!\!/H$ be the categorical
quotient with respect to $H$. The inclusion $S^G\subseteq S^H$
defines a morphism $\pi_{G/H}:Y/\!\!/H\to Y/\!\!/G$.
Obviously $\pi_{G/H}((Y/\!\!/H)^{[d]})\subseteq (Y/\!\!/G)^{[d]}$
because $\pi_{G/H}$ is $T$-equivariant (the inclusion $R^G\subseteq R^H$
respects the multi-grading).
We will prove that we have equality. Suppose that
$x\in (Y/\!\!/G)^{[d]}$. There exists $y\in\pi_G^{-1}(x)$ with a closed orbit,
and therefore there must exist a $g\in G$
such that $g\cdot y=\zeta\cdot y$ where $\zeta$ generates
a dense subgroup of $T^{[d]}$.
We have a multiplicative Jordan decomposition $g=g_sg_u$ where
$g_s$ is semi-simple and $g_u$ is unipotent.
Let $U$ be the Zariski closure of the group generated by $g_u$.
Because $U$ also is contained in the Zariski closure of the group generated
by $g$, we have $U\cdot y\subseteq \overline{T^{[d]}\cdot y}\subseteq 
\overline{T\cdot y}$. All elements in the $U$-orbit $U\cdot y$ have
the same stabilizers in $T$, and this shows that actually
$U\cdot y\subseteq T\cdot y$. Now $U$ is either the trivial group 
or the additive group. Since there is no way to embed the
affine line into a torus, the orbit $U\cdot y$ must be a point and
$g_u\cdot y =y$.
It follows that $g_s\cdot y=\zeta\cdot y$.
Now $g_s$ is conjugate to an element in the maximal torus of $G$.
Therefore, there exists  $a\in G$ such that $h:=ag_sa^{-1}\in H$.
We get that $h\cdot (a\cdot y)=\zeta\cdot (a\cdot y)$
and we get $z:=\pi_H(a\cdot y)\in (Y/\!\!/H)^{[d]}$.
We also have $\pi_{G/H}(z)=\pi_G(a\cdot y)=\pi_G(y)=x$.
In particular we now have that 
$$\dim (Y/\!\!/H)^{[d]}\geq\dim (Y/\!\!/G)^{[d]}.
$$
The theorem follows from Theorem~\ref{theo10}(a).
\end{proof}
\begin{remark}
We sketch an alternative  proof of Theorem~\ref{theo100} 
in the case where $H$ is the maximal torus.
Suppose that $I$ is an ideal of $S^G$. Then $I=J^G$ where
$J=IS$. Now from Weyl's character formula follows
$$
H(I,t)=\sum_{w\in W}(-1)^{\ell(w)}H((J\otimes \chi_{w(\rho)-\rho})^H,t)
$$
where $W$ is the Weyl group, $\ell(w)$ is the length of
an element $w\in W$, $\rho$ is the sum of the fundamental
weights and $\chi_{\lambda}$ is the character of $H$ corresponding
to the weight $\lambda$. Since $J$ is a finitely generated
$S$-module, $(J\otimes \chi_{w(\rho)-\rho})^H$ is a finitely
generated $S^H$-module. In particular,
the denominator of $H((J\otimes \chi_{w(\rho)-\rho})^H,t)$
divides the universal denominator of $S^H$.
It follows that the denominator of $H(I,t)$ divides
the universal denominator of $S^H$. Since
$I$ can be an arbitrary ideal, the universal
denominator of $S^G$ divides the universal denominator
of $S^H$.
\end{remark}
Suppose that $G$ is a linearly reductive group acting
regularly on a ring $S$ and $V$ is a irreducible representation
of $G$. Then $(S\otimes V)^G$ is a finitely generated $S^G$-module
and it is called a {\it module of covariants}. The universal
denominator of an invariant ring has the following interesting
interpretation in terms of modules of covariants.
\begin{theorem}
Suppose that $G$ is a linearly reductive algebraic group
acting homogeneously  and rationally on a multi-graded
polynomial ring $S$.
The universal denominator of $R=S^G$ is
$$
\lcm_V\{\denom((S\otimes V)^G,t)\},
$$
where $V$ runs over all irreducible representations of $G$.
\end{theorem}
\begin{proof}
Clearly $\denom(S\otimes V)^G,t)$ divides
$\udenom(R,t)$ because the module $(S\otimes V)^G$ 
of covariants is a finitely generated
$R$-module.
For any multi-graded $R$-module $M$, 
there exists a finite free $G$-equivariant graded minimal resolution
of $S$-modules
$$
0\to S\otimes V_k\to S\otimes V_{k-1}\to \cdots\to S\otimes V_0\to S\otimes_R M\to 0.
$$
where $V_i=\operatorname{Tor}^i(K,S\otimes_RM)$ is a multi-graded 
representation of $G$ for all $i$.
Taking $G$-invariants (which is an exact functor, since $G$
is linearly reductive) we get a free resolution
$$
0\to (S\otimes V_k)^G\to (S\otimes V_{k-1})^G\to \cdots\to 
(S\otimes V_0)^G\to (S\otimes_R M)^G=M\to 0.
$$
It follows
that
$$
H(M,t)=\sum_{i=0}^k(-1)^iH((S\otimes V_i)^G,t).
$$
This shows that $\udenom(R,t)$ divides 
$$
\lcm_i\{\denom(S\otimes V_i)^G,t)\},
$$
and the theorem follows.
\end{proof}
\section{Universal Denominators for Finite Groups Invariants}
Suppose now that $G$ is a finite group. We use the same notation
as in the previous section.
\begin{corollary}\label{cor102}
Suppose that $d$ is not divisible by the characteristic of $K$,
and $\zeta\in T^{[d]}$ generates a Zariski dense subgroup. Then we have
$$\pi^{-1}(X^{[d]})=\bigcup_{g\in G}Y^g_{\zeta}$$
\end{corollary}
\begin{proof}
This follows from Proposition~\ref{prop101} because every $G$-orbit is closed.
\end{proof}
\begin{theorem}
Suppose that the universal denominator of $H(R,t)$ is
$$
\prod_{d\in \N^r}\phi_d(t)^{m_d}.
$$
If the characteristic of $K$ does not divide $d$, then
$$
m_d=\max\{\dim Y^g_{\zeta}\mid g\in G\}.
$$
\end{theorem}
\begin{proof}
By Corollary~\ref{cor102} we have that 
$$
\pi^{-1}(X^{[d]})=\bigcup_{g\in G}Y^g_\zeta.
$$
Since $\pi$ is finite, 
$$m_d=\dim X^{[d]}=\dim \pi^{-1}(X^{[d]})=\max\{\dim Y^g_\zeta\mid g\in G\}.
$$
\end{proof}
\begin{corollary}\label{cor1020}
Assume that the characteristic of the base field $K$ is 0.
Suppose that $Y$ is a vector space on which the
finite group $G$ acts linearly.
Now $S=K[Y]$ is graded $S=\bigoplus_{d=0}^\infty S_d$.
For $g\in G$, let $\det(\id-tg)$ be the determinant of the endomorphism
$\id-tg$ where $t$ is an indeterminate. Then 
the universal denominator of $H(S^G,t)$ is
$$
\lcm\{\det(\id-tg)\mid g\in G\}
$$
where $\lcm$ is the least common multiple.
\end{corollary}
\begin{proof}
Let $\zeta$ be a $d$-th root of unity.
The multiplicity of $(1-\zeta t)$ in $\det(\id-tg)$
is exactly $\dim V^g_\zeta$.
It follows that the multiplicity of $(1-\zeta t)$ in
$$
\lcm\{\det(\id-tg)\mid g\in G\}
$$ 
is exactly 
$$
m_d=\max\{\dim V^g_\zeta\mid g\in G\}.
$$
\end{proof}
\begin{remark}
Corollary~\ref{cor1020} extends to the case that the characteristic
of the field $K$ does not divide the group order. In this more
general case $\det(\id-tg)$ has coefficients in $K$, not in $\Z$,
and $\Z$ is not necessarily a subring of $K$.
Since $\det(\id-tg)$ is a product of cyclotomic polynomials,
we can lift each cyclotomic polynomial with coefficients
in $K$ to a cyclotomic polynomial
with coefficients in $\Z$. In this way we can lift $\det(\id-tg)$
to a polynomial with coefficients in $\Z$, and Corollary~\ref{cor1020}
will make sense.
\end{remark}
\begin{remark}
Let us recall Molien's formula (see for example \cite[\S3.2]{DerksenKemper})
$$
H(S^G,t)=\sum_{g\in G}\frac{1}{\det(\id-tg)}.
$$
Molien's formula also holds when the characteristic of $K$ is
positive, but does
not divide the order of $G$. As in the previous remark, 
$\det(\id-tg)$ should be lifted to a polynomial 
with integer coefficients in that case.
It follows from Molien's formula that 
the denominator of $H(S^G,t)$
divides 
$$
\lcm\{\det(\id-tg)\mid g\in G\}.
$$
Our result here is stronger: the {\it universal\/} denominator
is always equal to this expression. However, it may happen
that the denominator of $H(S^G,t)$ is smaller than
the universal denominator of $H(S^G,t)$ (see Example~\ref{ex24}).
\end{remark}
\begin{example}
Consider the action of the symmetric group $\Sigma_n$ on $Y:=K^n$ where
$K$ is an algebraically closed field of characteristic $0$.
The coordinate ring $S:=K[Y]$ is graded: $S=\bigoplus_{d=0}^\infty S_d$.
Suppose $g\in \Sigma_n$ has cycle structure $(k_1,k_2,\dots,k_r)$ with $k_1\geq
k_2\geq\cdots\geq k_r\geq 1$ and $k_1+k_2+\cdots+k_r=n$, i.e.,
$k_1,k_2,\dots,k_r$ are the lengths of the cycles of the permutation $g$.
If $\zeta$ is a $d$-th root of unity, then
$$
\dim V^g_\zeta=\#\{i\mid \ d\mbox{ divides } k_i\}.
$$
We have
$$
\max_{g\in G}\{\dim V^{g}_{\zeta}\}=\left\lfloor \frac{n}{d}\right\rfloor.
$$
The maximum is reached if $g$ has $\lfloor {n/d}\rfloor$ $d$-cycles.
So the universal denominator of $S^{\Sigma_n}$ is
$$
\prod_{d=0}^\infty \phi_{d}(t)^{\lfloor n/d\rfloor}=
(1-t)(1-t^2)\cdots (1-t^n).
$$

We can directly verify this.
It is well known that $S^{\Sigma_n}=K[e_1,e_2,\dots,e_n]$ where $e_i$ is
the $i$-th elementary symmetric function of degree $i$.
The Hilbert series of $S^{\Sigma_n}$ is
$$
H(S^{\Sigma_n},t)=\frac{1}{(1-t)(1-t^2)\cdots (1-t^n)}.
$$
Since $S^{\Sigma_n}$ is a polynomial ring, we also have that
the denominator of any finitely generated $S^G$-module has
a denominator dividing $(1-t)(1-t^2)\cdots (1-t^n)$. 
So indeed,  $(1-t)(1-t^2)\cdots (1-t^n)$ is
the universal denominator of $S^{\Sigma_n}$.
\end{example}
\begin{example}\label{ex24}
Let $A_n\subset \Sigma_n$ be the alternating group
acting on $K^n$ as in the previous example.
We have 
$$
\udenom(S^{A_n},t)=\prod_{d\in \N}\phi_d(t)^{m_d}
$$
where
$$
m_d:=\max\{\dim V^g_\zeta\mid g\in A_n\}
$$
for any primitive $d$-th root of unity $\zeta$.
We now have that
$$
m_d=\left\{\begin{array}{ll}
\lfloor n/d\rfloor & \mbox{if $d$ is odd, $\lfloor n/d\rfloor$
is even or $d\lfloor n/d\rfloor\leq n-2$,}\\
\lfloor n/d\rfloor -1 & \mbox{otherwise.}\\
\end{array}\right. 
$$
If $d$ is odd or $\lfloor n/d\rfloor$ is even, then this maximum is
reached by taking for $g$ a product of $\lfloor n/d\rfloor$
disjoint $d$-cycles (this indeed gives an even permutation).
Suppose now that $d$ is even and $\lfloor n/d\rfloor$ is odd.
If $d\lfloor n/d\rfloor\leq n-2$, then the maximum is reached
by taking for $g$ a product
of $\lfloor n/d\rfloor$ $d$-cycles and a $2$-cycle. In
any other case, the maximum is achieved by taking $\lfloor n/d\rfloor-1$
$d$-cycles.

Let us compute the denominator of $H(S^{A_n},t)$.
We have 
$$
S^{A_n}=S^{\Sigma_n}\oplus S^{\Sigma_n}\Delta
$$
where 
$$
\Delta=\prod_{1\leq i<j\leq n}(x_i-x_j)
$$
is an $A_n$-invariant of degree $n(n-1)/2$.
In particular,
$$
H(S^{A_n},t)=\frac{1+t^{n(n-1)/2}}{(1-t)(1-t^2)\cdots (1-t^n)}.
$$
Note that
$$
1+t^{n(n-1)/2}=\frac{1-t^{n(n-1)}}{1-t^{n(n-1)/2}}=\prod_{d} 
\phi_d(t)
$$
where the product is over all $d$ with $d\mid n(n-1)$ and $d\not\,\mid n(n-1)/2$.
These are exactly all $d$ such that $n(n-1)/d$ is an odd positive integer.
So we have
$$
\denom(S^{A_n},t)=\prod_{d\in \N}\phi_{d}(t)^{a_d}
$$
where
$$
a_d=\left\{\begin{array}{ll}
\max\{0,\lfloor n/d\rfloor-1\} & \mbox{if $n(n-1)/d$ is an odd positive integer,}\\
\lfloor n/d\rfloor & \mbox{otherwise.}\end{array}\right.
$$
The reader may check that $a_d\leq m_d$ for all $d$.
Note that in this example, the denominator is not always equal
to the universal denominator. For example, if we take $n=10$, we have
$$
\denom(R,t)=\phi_1(t)^{10}\phi_2(t)^4\phi_3(t)^3\phi_4(t)^2
\phi_5(t)^2\phi_7(t)\phi_8(t)\phi_9(t).
$$
and
$$
\udenom(R,t)=\phi_1(t)^{10}\phi_2(t)^4\phi_3(t)^3\phi_4(t)^2
\phi_5(t)^2\phi_6(t)\phi_7(t)
\phi_8(t)\phi_9(t).
$$
Note that the universal denominator has an additional factor
$\phi_6(t)=1-t+t^2$. This means that although $\phi_6(t)$ does not
appear in the denominator of $H(R,t)$, it does
appear in the denominator
$H(M,t)$ for some finitely generated $R$-module.
We will describe such a module. Define the ideals
$$
{\mathfrak q}=(e_1,e_2,e_3,e_4,e_5,e_7,e_8,e_9,e_{10})
$$
and
$$
{\mathfrak p}={\mathfrak q}+(\Delta)
$$
of $R$. 
We claim that $R/{\mathfrak p}$ is the polynomial ring 
generated by (the image of) $e_6$.

Note that $e_6$ does not lie in the radical
ideal $\sqrt{S{\mathfrak q}}$ of $S{\mathfrak q}$, 
since 
$$e_1=e_2=\cdots=e_5=e_7=\cdots=e_{10}=0$$
does not imply that $e_6=0$.
Also note that $\Delta$ does lie in  $\sqrt{S{\mathfrak q}}$
(since the polynomial $X^{10}+e_6X^4$ has a multiple zero, namely
0). 
Obviously, $R/{\mathfrak p}$ is generated by $e_6$. Also,
no power of $e_6$ lies in ${\mathfrak p}$
because $e_6\not\in \sqrt{S{\mathfrak p}}=\sqrt{S{\mathfrak q}}$. It follows
that $R/{\mathfrak p}$ is the polynomial ring generated by
the invariant $e_6$ of degree $6$.
In particular, we have
$$
H(R/{\mathfrak p},t)=\frac{1}{1-t^6}.
$$
so $\phi_6(t)$ divides the denominator of the Hilbert series
of the $R$-module $R/{\mathfrak p}$.
\end{example}
\section{Universal Denominators for Invariants of  Tori}
Let $K$ be an algebraically closed base field of characteristic 0.
Let $G:=(K^\star)^l$ be an $l$-dimensional torus
acting 
on $Y:=K^n$ diagonally. The coordinate ring
of $Y$ is $S:=K[Y]=K[y_1,\dots,y_n]$.
If $g=(g_1,\dots,g_l)\in G$
and $\omega=(\omega_1,\dots,\omega_l)\in \Z^l$ then we write
$g^{\omega}=g_1^{\omega_1}\cdots g_l^{\omega_l}$.
Since $G$ acts diagonally, there exist $\omega_1,\omega_2,\dots,\omega_n\in
\Z^l$ such that $g\cdot y_i=g^{\omega_i}\cdot y_i$.
We assume that $K[Y]$ is multi-graded such that $y_i$ is homogeneous of 
degree $d_i$ for all $i$ with $d_1,\dots,d_n\in \N^r$.
Let $R=S^G$ be the invariant ring. In this section, we would like to describe
the universal denominator of the Hilbert series $H(R,t)$.

Suppose that $I\subseteq \{1,2,\dots,n\}$ is a subset. We define
$\Omega_I$ to be the $\Z$-module generated by $\omega_i$, $i\in I$.
We let $\cx_I$ to be the polytope spanned by $\omega_i$, $i\in I$
inside $\Omega_I\otimes_\Z\R$. We let $M_I\subseteq \Z^l\times \Z^r$
be the module generated by all $(\omega_i,d_i)$, $i\in I$.
\begin{theorem}
Write 
$$
\udenom(R,t)=\prod_{d\in \N^r}\phi_d(t)^{m_d}.
$$
We have
$$
m_d=\max_I\{\#I-\rank \Omega_I\}
$$
where $I$ runs over all subsets of $\{1,2,\dots,n\}$ such 
that $0$ lies in the interior of $\cx_I$ and
$$
M_I\cap \{0\}\times \Z^r\subseteq \{0\}\times \Z d.
$$
\end{theorem}
\begin{proof}
Let us define
$$
Y_I=\{(y_1,y_2,\dots,y_n)\in Y\mid y_i=0 \mbox{ for all $i\not\in I$}\}.
$$
Let $Z^{[d]}\subseteq Y$ be the closure of the (constructible) set of all 
closed orbits in 
$\pi^{-1}(X^{[d]})$.
We claim that $Z^{[d]}$ is a union of sets of the form
 $Y_I$ with $I\subseteq\{1,2,\dots,n\}$.
 

Take $\zeta$ a generator of a dense
Zariski dense subgroup of $T^{[d]}$.
Note that $Z^{[d]}$ is contained in 
the set of closed orbits in $\bigcup_{g\in S}Y^g_\zeta$
where $S$ ranges over a subset of $G$ by Proposition~\ref{prop101}. 
Each $Y^g_\zeta$
is of the form $Y_I$ for some $I\subseteq \{1,2,\dots,n\}$
and the closure of the union of closed orbits in $Y_I$
is equal to $Y_J$ for some subset $J\subseteq I$
(namely, take the largest face $F$ of $\cx_I$ such that
$0$ lies in the interior of $F$, then let $J$
be the set of all $i$ with $\omega_i\in F$).
This shows that $Z^{[d]}$ is a finite  union of sets of the
form $Y_I$ with $I\subseteq \{1,2,\dots,n\}$
such that the general orbit of $Y_I$ is closed.

Now it follows that 
$$m_d=\dim X^{[d]}=\dim \pi(Z^{[d]})$$ 
is the maximum of all $\dim \pi(Y_I)$ with $Y_I\subseteq Z^{[d]}$
and such that the general orbit of $Y_I$ is closed.
If we restrict $\pi$ to such a set $Y_I$, 
then for $y\in Y_I$ in general position
$\pi^{-1}\pi(y)$ is the (closed) orbit of $y$ which
has dimension $\rank \Omega_I$. Since $\dim Y_I=\#I$ this shows
that
$$
\dim \pi(Y_I)=\#I-\rank I.
$$

Note that the general orbit of $Y_I$ is closed if and only if
$0$ lies in the interior of $\cx_I$. Let us assume
that the general orbit in $Y_I$ is closed.
We claim that $Y_I\subseteq Z^{[d]}$ if and only if
$$
M_I\cap \{0\}\times \Z^r\subseteq \{0\}\times \Z d.
$$
First note that $Y_I\subseteq Z^{[d]}$ if and only if
$Y_I\subseteq Y^g_\zeta$ for some $g\in G$.
Now $Y_I\subseteq Y^g_\zeta$ if and only if
the system of equations
\begin{equation}\label{eqg}
g^{\omega_i}=\zeta^{d_i},\quad i\in I
\end{equation}
has a solution $g\in G$.
Let $A\subseteq (K^\star)^I$ be the subgroup
of all $(g^{\omega_i},i\in I)$, with $g\in G$.
Let $B\subseteq (K^\star)^I$ be the subgroup
of all $(t^{d_i},i\in I)$, with $t\in T^{[d]}$.
Now the system of equations~(\ref{eqg}) has a solution
if and only if $A\supseteq B$.
The coordinate ring of the algebraic group $(K^\star)^I$
is 
$$
R=K[z_i,z_i^{-1}, i\in I]
$$
The vanishing ideals $I(A),I(B)\subset R$ are generated by elements
of the form $m-1$ where $m$ is a Laurent monomial in $z_i,i\in I$.
Now $A\supseteq B$ if and only if $I(A)\subseteq I(B)$.
Suppose that
$$
m=\prod_{i\in I}z_i^{a_i}
$$
then $m-1\in I(A)$ if and only if 
$$
\prod_{i\in I} g^{a_i\omega_i}=1
$$
for all $g\in G$ which is equivalent to $\sum_{i\in I} a_i\omega_i=0$.
On the other hand
$m-1\in I(B)$ if and only if
$$
\prod_{i\in I} t^{a_id_i}
$$
lies in the ideal generated by $t^d-1$.
So $m-1\in I(B)$ if and only if $\sum_{i\in I}a_id_i$ is 
a multiple of $d$.
This shows that $I(A)\subseteq I(B)$ if and only if
$$
M_I\cap \{0\}\times \Z^r\subseteq \{0\}\times \Z d.
$$
\end{proof}

Let us now consider the case that $r=1$, and that all
variables $y_1,y_2,\dots,y_n$ have degree 1.
Suppose that $\Theta\subseteq \Z^r$ is a $\Z$-submodule 
and $a\in \Z^r$.
Then we define 
$$
I_{a+\Theta}=\{i\in \{1,2,\dots,n\}\mid \omega_i\in a+\Theta\}.
$$
\begin{corollary}
Let us write
$$
\udenom(R,t)=\prod_{d\in \N}\phi_d(t)^{m_d},
$$
then we have
$$
m_d=\max_{\lambda,\Theta}\{\#I_{\lambda+\Theta}-\rank\Theta\}
$$
where
$\Theta$ runs over all $\Z$-submodules of $\Z^r$ 
and $\lambda$ runs over all torsion elements of $\Z^r/\Theta$ for which
$$
(\Z \lambda+\Theta)/\Theta\cong \Z/d\Z
$$
and
$0$ lies in the interior of $\cx_{I_{\lambda+\Theta}}$.
\end{corollary}
\begin{proof}
If we apply the previous theorem, we will see that we
only need to consider subset $I$ of the form $I_{\lambda+\Theta}$.
Suppose that $I=I_{\lambda+\Theta}$ and
$$
(\Z \lambda+\Theta)/\Theta\cong \Z/d\Z.
$$
If $\sum_{i\in I} a_i \omega_i=0$,
then $(\sum_{i\in I}a_i)\lambda\in \Theta$
and by the above isomorphism $d$ must divide $\sum_{i\in I}a_i$.
We have proven that
$$
M_I\cap \{0\}\times \Z\subseteq \Z d.
$$

Conversely, suppose that for some subset $I\subseteq \{1,2,\dots,n\}$,
 $0$ lies in the interior of $\cx_I$ and
$$
M_I\cap \{0\}\times \Z\subseteq \Z d.
$$
 Put $\lambda=\omega_i$ for some $i\in I$.
 Let $\Theta$ be the $\Z$-module generated
 by all $\omega_i-\omega_j$, $i,j\in I$ and $d\lambda$.
 Clearly we have
 $$
 (\Z\lambda+\Theta)/\Theta\cong\Z/d\Z.
 $$
 Now $I\subseteq I_{\lambda+\Theta}$, $0$ lies in the relative
 interior of $I_{\lambda+\Theta}$ and
 $$
 M_{\lambda+\Theta}\cap \{0\}\times \Z\subseteq \Z d.
 $$
Also, obviously
$$
\#I-\rank \Omega_I=\#I-\rank \Theta\leq
\#I_{\lambda+\Theta}-\rank \Theta.
$$
\end{proof}

Let us now again specialize. Let us assume that $l=1$, i.e., 
$G=K^\star$ is the one-dimensional torus.
\begin{corollary}
Write 
$$
\udenom(R,t)=\prod_{d\in \N}\phi_d(t)^{m_d}.
$$
If $\omega_i\geq 0$ for all $i\in \{1,2,\dots,n\}$ or $\omega_i\leq 0$ 
for all $i\in \{1,2,\dots,n\}$
then
$m_1=\{i\mid \omega_i=0\}$. Otherwise $m_1=n-1$.

For $d\geq 2$ we have
$$
m_d=\max_{a,b}\{\#\{i\mid \omega_i\in ab+ad\Z\}\}-1
$$
where $a$ and $b$ run over all integers such that $b$ is relatively
prime to $d$ and
$$
(ab+ad\Z)\cap\{\omega_1,\dots,\omega_n\}
$$
contains at least one positive and one negative integer.
(If we assume that $\omega_1\leq \omega_2\leq\cdots\leq \omega_n$
then one only needs those $a,b$ such that $ad\leq \omega_n-\omega_1$
and $1\leq b\leq d-1$ such that $b$ and $d$ are relatively prime).

\end{corollary}
\begin{proof}
We apply the previous corollary. For $d=1$, the rank of $\Theta$
may be $0$ or $1$. If the rank of $\Theta$ is 0 then
$\lambda$ must be 0. Then $\# I_{\lambda+\Theta}-\rank \Theta$ is the number
of zero weights. 
If $\Theta$ has rank 1, then the maximum is achieved if
$\Theta=\Z$. The condition that 0 is in the interior
of $\cx_I$ means that there are positive and negative weights.
In this case we get $\#I_{\lambda+\Theta}-\rank \Theta=n-1$.

Suppose that $d>1$. Now $\Theta$ can only have rank 1.
Again the condition that 0 is in the  interior
of $I_{\lambda+\Theta}$ means that $\lambda+\Theta$
contains both positive and negative weights.
Write $\Z\lambda+\Theta=a\Z$ for some nonnegative $a\in \Z$.
Then we have $\Theta=ad\Z$ and $\lambda=ab$ with $b$ relatively prime to $d$.
\end{proof}

\begin{example}
Let $G=K^\star$ act on $K^5$ with weights $-3,-2,2,5,6$.
Since there are positive and negative weights,
we get $m_1=5-1=4$.
We compute
$$
\{\#\{i\mid \omega_i\in ab+ad\Z\}\}
$$
for all $a,b,d$ with $1\leq ad\leq 9$ and $1\leq b\leq d-1$
with $b,d$ relative prime and such
that
$$(ab+ad\Z)\cap\{-3,-2,2,5,6\}
$$
contains negative and positive integers.
$$
\begin{matrix}
(d=2) & 1+2\Z:2 & 2+4\Z:3 & 3+6\Z:2 & \Rightarrow m_2=2\\
(d=3) & 6+9\Z:2 & & & \Rightarrow m_3=1\\
(d=4) & 1+4\Z:2 & 2+8\Z:2 & & \Rightarrow m_4=1\\
(d=5) & 2+5\Z:2 & & & \Rightarrow m_5=1\\
(d=6) & - & & & \Rightarrow m_6=0\\
(d=7) & 5+7\Z:2 & & & \Rightarrow m_7=1\\
(d=8) & 5+8\Z:2 & & & \Rightarrow m_8=1
\end{matrix}
$$
The universal denominator of $H(K[Y]^G,t)$ is therefore equal to
$$
\phi_1(t)^4\phi_2(t)^2\phi_3(t)\phi_4(t)\phi_5(t)\phi_7(t)\phi_8(t).
$$
\end{example}
\begin{example}
Consider the action of the multiplicative group $K^\star$
on $Y_n:=K^{n+1}$ with weights
$$
-n,-n+2,\dots,n-2,n.
$$
We will describe the universal denominator. 

Let us first assume that {\bf $n$ is odd}. Put $S_n=\{-n,-n+2,\dots,n\}$.
Note that $m_1=(n+1)-1=n$.

Choose $d\geq 2$.
If the cardinality of
$(ab+ad)\Z\cap S_n$ is at most 1, 
then it cannot contain both a positive and a negative element.
If the cardinality of
$(ab+ad)\Z\cap S_n$ is $\geq 2$, 
then this intersection will contain automatically a positive and
a negative element. If $a\geq 2$, then
$$
(ab+ad\Z)\cap S_n\subseteq a\big(
 (b+d\Z) \cap S_n\big)
$$
because $S_n$ is closed under taking divisors.
It is therefore clear that we only need to consider the case $a=1$.
We now have  to maximalize
$$
\# (b+d\Z)\cap S_n
$$
over all $b$ relatively prime to $d$.

Let us assume that {\bf $d$ is odd}.
If $d$ divides $n$, then $b$ cannot be divisible by $d$.
It follows that
$$
\max_b\# (b+d\Z)\cap S_n=n/d.
$$
If $d$ does not divide $n$, then
$$
\# (b+d\Z)\cap S_n\leq \lceil n/d\rceil.
$$
Write $n=kd+e$ with $0<e<d$, then $k+1=\lceil n/d\rceil$. We have equality
if $b=n-e+2$ if $e$ is even, and $b=n-e+1$ if $e$ is odd.
We have proven that
$$
m_d=\lceil n/d\rceil-1,\quad \mbox{$d$ odd, $d\geq 3$}.
$$

Let us assume that {\bf $d$ is even}.
If $d=2$ then 
$$
\# (1+2\Z)\cap S_n=n+1
$$
and $m_2=n$.
Assume $d\geq 4$. If $d/2$ divides $n$ then $d/2$ cannot
divide $b$. It follows that in this situation
we have
$$
\# (b+d\Z)\cap S_n=2n/d,
$$
Otherwise, we have that
$$
\# (b+d\Z)\cap S_n\leq \lceil 2n/d\rceil.
$$
If we write $n=k(d/2)+e$ with $0<e<(d/2)$
then we have equality for $b=n-e+1$ if $e$ is odd and
for $b=n-e+2$ if $e$ is even.
It follows that
$$
m_d=\lceil 2n/d\rceil-1,\quad\mbox{$d$ even, $d\geq 4$.}
$$

Let us now assume that {\bf $n$ is even}.
Again we have $m_1=n$. 

Choose $d\geq 2$. 
Write $a=pq$ where $q=1$ if $a$ is odd and $q=2$ if $a$ is even.
Then it is easy to check that
$$
(ab+ad\Z)\cap S_n\subseteq p\big((qb+qd\Z)\cap S_n\big).
$$
Therefore we only need to consider the cases $a=1$ and $a=2$.

Assume that {\bf $d$ is even}. In that case $b$ must be odd
and therefore we must have $a=2$ (otherwise the intersection is empty).
One can easily check that we get $m_2=n/2$ if $n\equiv 2\bmod 4$ 
and $m_2=n/2-1$ if $n\equiv 0\bmod 4$, because $b$ is odd.
Suppose that $d\geq 4$.
If $d$ divides $n$ then $m_d=n/d-1$ (by similar reasoning as before).
Otherwise, write $n/2=k(d/2)+e$ with $0<e<d/2$.
We have
$$
\# 2(b+d\Z)\cap S_n\leq \lceil n/d\rceil.
$$
If $n$ is not divisible by $4$,
then we have equality for $b=n/2-e+1$ if $e$ is odd,
or $b=n/2-e+2$ if $e$ is even. If $e\neq 1$, i.e., $n-2$ is
not divisible by $d$, then we have equality for $b=n/2-e+2$.

In the remaining case, $n$ divisible by $4$ and $d$ divides $n-2$,
one can easily see that we cannot have equality. So if $n$
and $d$ are both even then, 

$$
m_d=\left\{\begin{array}{ll}
\lceil n/d\rceil-2 & \mbox{if $4\mid n$, $d$ even and $d\mid n-2$;}\\
\lceil n/d\rceil-1 & \mbox{otherwise.}\\
\end{array}\right.
$$

Assume now that {\bf $d$ is odd}. Since
$$
(b+d\Z)\cap S_n\subseteq (2(b/2)+2d\Z)\cap S_n
$$
if $b$ is even and
$$
(b+d\Z)\cap S_n\subseteq (2((b+d)/2)+2d\Z)\cap S_n
$$
if $b$ is odd, we may only consider the case that $a=2$.
Similar reasoning as before gives us
$$
m_d=\lceil n/d\rceil-1.
$$

The universal denominator of $K[Y_n]^{K^\star}$ is equal to
$$
\left\{\begin{array}{ll}
(1-t^2)^2(1-t^4)\cdots (1-t^{2n-2}) & \mbox{$n$ odd;}\\
(1-t)(1-t^2)^2(1-t^3)(1-t^4)\cdots (1-t^{n-1}) & \mbox{$n\equiv 2\bmod 4$};\\
(1-t)(1-t^2)^2(1-t^3)(1-t^4)\cdots (1-t^{n-3})(1-t^{n/2-1})(1-t^{n-1})
& \mbox{$n \equiv 0\bmod 4$}.\\
\end{array}\right.
$$

\end{example}
\section{Universal Denominators  for Binary Forms}
Let $G:=\SL_2$ act on $Y$ where 
$$
Y=\{a_0x^n+a_1x^{n-1}y+\cdots+a_ny^n\}
$$
is the set of binary forms of degree $n$ ($n\geq 3$).
Let $S=K[a_0,a_1,\dots,a_n]$ and let $R=S^{\SL_2}$.
Let $H$ be the maximal 1-dimensional torus of $\SL_2$
consisting of all diagonal matrices:
$$
H=\left\{\begin{pmatrix}\lambda & 0 \\ 0  & \lambda^{-1}\end{pmatrix},
\lambda\in K^\star\right\}.
$$
We already computed the universal denominator 
$$\udenom(S^{H},t)=\prod_{d\in \N}\phi_d(t)^{m_d}$$ 
Let us write
$$
\udenom(S^{\SL_2},t)=\prod_{d\in \N}\phi_d(t)^{u_d}.
$$
We have $u_1=n-1$ if $n$ is odd  and $u_1=u_2=n-1$ if $n$ is even.
Suppose that $d\geq 3$ and $n$ is even or $d\geq 2$ and $n$ is odd.
We will prove that $u_d=m_d$.
We already know
that $u_d\leq m_d$ by Theorem~\ref{theo100}. Let us prove that $u_d\geq m_d$.
If $m_d=0$ then there is nothing to prove. We will assume
that $m_d\geq 1$.

Let $\zeta$ be a $d$-th primitive root of unity.
There exists an 
$$h=\begin{pmatrix}\lambda & 0 \\ 0  & \lambda^{-1}\end{pmatrix}\in H$$ 
such that
$$
\dim Y^h_\zeta\geq m_d+1\geq 2.
$$
The space $Y^h_{\zeta}$ is the set of all polynomials
$$
f=a_0x^n+a_1x^{n-1}y+\cdots+a_ny^n
$$
for which $a_i=0$ whenever $\lambda^{n-2i}\neq \zeta$.
We see that $\lambda$ is an $(pd)$-th root of unity for some positive
integer $p$. A general polynomial $f\in Y^h_\zeta$ has at least
3 distinct roots and no root with multiplicity $\geq n/2$.
It follows that the $G$-orbit of such a general polynomial $f$
is closed and 3-dimensional. The dimension of a general
fiber $\pi_G^{-1}\pi_G(f)$ is 3 and
$$
\dim \overline{\pi_G(G\cdot Y_{\zeta}^h)}=\dim \overline{G\cdot Y_{\zeta}^h}-3.
$$
Consider the morphism 
$$
\phi:G\times Y_{\zeta}^h\to Y
$$
defined by
$$
(g,f)\mapsto g\cdot f.
$$
Again, let $f\in Y_{\zeta}^h$ be a polynomial in general position.
One can easily verify that the rank of the differential $d\phi$
at $(e,f)$ has rank $\dim Y_{\zeta}^h+2\geq m_d+3$.
We have $\dim \overline{G\cdot Y_{\zeta}^h}\geq m_d+3$.
Since
$\pi_G(\overline{G\cdot Y_{\zeta}^h})\subseteq (Y/\!\!/G)^{[d]}$
we obtain
$$
u_d=\dim (Y/\!\!/G)^{[d]}\geq \dim \overline{\pi_G(G\cdot Y_{\zeta}^h})
=\dim \overline{G\cdot Y_{\zeta}^h}-3\geq m_d
$$
 and therefore $u_d=m_d$.
 Finally we obtain that the universal denominator of $K[Y_n]^{\SL_2}$
 is equal to
 $$
 \left\{\begin{array}{ll}
 (1-t^4)(1-t^6)\cdots (1-t^{2n-2}) & \mbox{if $n$ is odd;}\\
 (1+t)(1-t^2)(1-t^3)\cdots (1-t^{n-1}) & \mbox{if $n\equiv 2 \bmod 4$;}\\
 (1+t)(1-t^2)(1-t^3)\cdots (1-t^{n-3})(1-t^{n/2-1})(1-t^{n-1})
 & \mbox{if $n\equiv 0\bmod 4$.}\\
 \end{array}\right.
 $$
 These formulas for the universal denominators are the same
 as the formulas for the denominators of the Hilbert series of
 binary forms as conjectured by Dixmier (see~\cite{Di}).
 It could happen, of course, that the universal denominator is not
 equal to the denominator due to some unfortunate accidental cancelling
in the Hilbert series.
The universal denominator is more interesting
than the denominator of the Hilbert series, since it has
a geometric interpretation and nice properties.
 Therefore, although Dixmier's conjecture is still unsettled,
 the formula for the universal denominator here may is just
 as satisfying as a positive answer to Dixmier's conjecture
 for the denominator of the Hilbert series of binary forms.

\end{document}